% \magnification=1200

% qualche utile definizione

\input epsf.tex

\def\2{{1\over 2}}

\def\d{\delta}

\def\b{\beta}
\def\g{\gamma}

\def\e{\epsilon}

\def\D{\Delta}

\def\fun#1#2#3{#1\colon #2\rightarrow #3}

\def\abs#1{\vert #1 \vert}
\def\frac#1#2{{{#1} \over {#2}}}

\def\st{\;\colon\;}
\def\tends{\rightarrow}

\def\dr{ {\rm d} }

\def\R{{\bf R}}

\def\Z{{\bf Z}}

\def\T{{\bf T}}

\def\thm#1{\vskip 1 pc\noindent{\bf Theorem #1.\quad}\sl}
\def\lem#1{\vskip 1 pc\noindent{\bf Lemma #1.\quad}\sl}
\def\prop#1{\vskip 1 pc\noindent{\bf Proposition #1.\quad}\sl}
\def\cor#1{\vskip 1 pc\noindent{\bf Corollary #1.\quad}\sl}

\def\proof{\rm\vskip 1 pc\noindent{\bf Proof.\quad}}
\def\fin{\par\hfill $\backslash\backslash\backslash$\vskip 1 pc}
\def\txt#1{\quad\hbox{#1}\quad}

\def\L{{\cal L}}

\def\G{\Gamma}

\def\r{\rho}

\def\cin#1{\2\abs{{#1}}^2}

\def\2{\frac{1}{2}}
\def\inn#1#2{{\langle #1 ,#2\rangle}}

\def\part{{\partial_{x}}}

\def\pprime{{{}^\prime{}^\prime}}
\def\ec{{\cal E}}

\def\fc{{\cal F}}

\def\wc{{\cal W}}
\def\dc{{\cal D}}
\def\pc{{\cal P}}
\def\tc{{\cal T}}

\def\ini{{(t,x)}}

%\def\exp{{\rm exp}}
%\def\freq{{\frac{\bar k}{\bar T}}}

% pagina di titolo 

\baselineskip= 17.2pt plus 0.6pt
% \hsize=14truecm
% \hoffset=2truecm
% \vsize=8truein
% \voffset=0.5truein
\font\titlefont=cmr17
\centerline{\titlefont The stochastic value function}
\vskip 1 pc
\centerline{\titlefont in metric measure spaces}
%\vskip 1 pc
% \footnote{}{\rm Supported by 
% Italian Ministry of Education}}
\vskip 4pc
\font\titlefont=cmr12
\centerline{         \titlefont {Ugo Bessi}\footnote*{{\rm 
Dipartimento di Matematica, Universit\`a\ Roma Tre, Largo S. 
Leonardo Murialdo, 00146 Roma, Italy.}}   }{}\footnote{}{
{{\tt email:} {\tt bessi@matrm3.mat.uniroma3.it} Work partially supported by the PRIN2009 grant "Critical Point Theory and Perturbative Methods for Nonlinear Differential Equations}} 
\vskip 0.5 pc
%\centerline{{\tt email:} {\tt bessi@matrm3.mat.uniroma3.it}}
 
\par
\vskip 2pc
\centerline{\bf Abstract}
Let $(S,d)$ be a compact metric space and let $m$ be a Borel probability measure on 
$(S,d)$. We shall prove that, if $(S,d,m)$ is a $RCD(K,\infty)$ space,  then the stochastic value function satisfies the viscous Hamilton-Jacobi equation, exactly as in Fleming's theorem on $\R^d$.

\vskip 2 pc
%\magnification=1200
\centerline{\bf  Introduction}
\vskip 1 pc

Let $\T^d\colon=\frac{\R^d}{\Z^d}$ denote the $d$-dimensional torus, let 
$w_0\in C^2(\T^d)$ and let $\fun{\fc}{(-\infty,0]\times\T^d}{\R}$ be a continuos potential. It is well-known ([11]) that the solution of the Schr\"odinger equation, backward in time, 
$$\left\{
\eqalign{
\partial_t w(t,x)&+\2\D w(t,x)+\fc(t,x)w(t,x)=0,\quad t\le 0\cr
w(0,x)&=w_0(x)
}
\right.  $$
is given by the Feynman-Kac formula
$$w(t,x)=E^{(t,x)}\left\{
e^{\int_t^0\fc(\tau, B^{(t,x)}_\tau)\dr\tau} w_0(B^{(t,x)}_0)
\right\}$$
where $B^{(t,x)}_\tau$ is the Brownian motion starting at $B^{(t,x)}_t\equiv x$ and 
$E^{(t,x)}$ denotes expectation with respect to Wiener's measure. Moreover, if 
$w(t,\cdot)>0$ for all $t\in(-\infty,0]$ (which happens if the final condition $w_0$ is positive), we can define
$$u(t,x)\colon=-\log w(t,x)  \eqno (1)$$
which solves the viscous Hamilton-Jacobi equation, backward in time, 
$$\left\{
\eqalign{
\partial_t u(t,x)&+\2\D u(t,x)-\cin{\nabla u(t,x)}-\fc(t,x)=0,\quad t\le 0\cr
u(0,x)&=u_0(x)\colon=-\log w_0(x)  .  
}
\right.   \eqno (2)$$
On the other side ([13]), $u$ is a value function, i. e. 
$$u(t,x)=\inf E^\ini\left\{
\int_t^0\left[
\cin{Y(\tau,x_\tau)}-\fc(\tau,x_\tau)
\right]   \dr\tau+u_0(x_0)
\right\}  , $$
where $x_\tau$ solves the Stochastic Differential Equation
$$\left\{
\eqalign{
\dr x_\tau&= Y(\tau,x_\tau)\dr\tau+\dr B_\tau^\ini\cr
x_t&=x,
}
\right.   \eqno (3)$$
and the $\inf$ is over all smooth vector fields $Y$. 

It is natural to ask whether some of these facts remain true in a more general setting. If we look at the various ingredients of Fleming's proof, the Brownian motion is the one with the longest history. Brownian motions on fractals have been studied since the Eighties (see [7], [8] and references therein); the crucial connection with Dirichlet forms was proposed in [17]. The "minimal" requirement to have a Brownian motion is the following: $(S,d)$ is a compact metric space, $m$ is a Borel probability measure on $S$, positive on open sets, and $\ec$ is a strongly local Dirichlet form on $L^2(S,m)$ (see section 1 below for the precise definitions). By [14], this implies the existence of a Brownian motion starting from $m$ a. e. $x\in S$.  

Next, we have to make sense of (3) or, equivalently, of the Fokker-Planck equation. We recall that, on $\T^d$, $\mu$ satisfies the weak version of the Fokker-Planck equation with drift $X$ on the interval $(t,0)$ if
$$\int_t^0\dr s\int_{\T^d}[
\partial_s\phi+\2\D\phi+\inn{X}{\nabla\phi}
]  \dr\mu_s  =0 $$
for all test functions $\phi\in C^\infty_0((t,0)\times\T^d)$. All of this translates to our setting: $\D$ becomes $\D_\ec$, the "Laplacian" associated with the Dirichlet form $\ec$. As for the internal product $\inn{X}{\nabla\phi}$ the theory of Dirichlet forms provides an object which behaves similarly: it is called the carr\'e\ de champs, and we shall suppose that the carr\'e\ de champs is defined on $\dc(\ec)$, the domain of $\ec$. As we shall see in section 3 below, one can also define a class of test function $\tc$, namely the functions $\phi$ such that $\phi\in C^1([t,0],L^2(S,m))\cap L^\infty([t,0],\dc(\D_\ec))$. 

This setting is sufficient to prove points 1) and 2) of theorem 1 below; if we want to go farther, we need to prove that the function $u$ defined in (1) above satisfies the Hamilton-Jacobi equation (2). In other words, we need information on the Laplacian 
$\D_\ec(\log w)$. It turns out that $\D_\ec(\log w(t,\cdot))\in L^2(S,m)$ if the carr\'e\ de champs of $w(t,\cdot)$ belongs to $L^2$. That's why we need our last ingredient, i. e. that $(S,d,m)$ is a $RCD(K,\infty)$ space and that $\ec$ is the double of Cheeger's energy: in this setting, it is standard that $w(t,\cdot)$ is Lipschitz and that the carr\'e\ de champs of 
$w(t,\cdot)$ is bounded by its Lipschitz constant. Using these facts and the method of [13], we shall be able to prove one inequality of formula (7) below. For the opposite inequality, we need to solve the Fokker-Planck equation with drift $\nabla u$; again, the fact that $u$ is Lipschitz will be essential. 

We shall use the strategy just outlined to prove the following theorem; we refer the reader to the next sections for the definitions of the various terms appearing in it. 

\thm{1}  Let $(S,d)$ be compact and let $m$ be a Borel probability measure on $S$, positive on open sets. Let us suppose that $(S,d,m)$ is a $RCD(K,\infty)$ space and let us call $\ec$ the natural Dirichlet form on $L^2(S,m)$, i. e. the double of Cheeger's energy. Let $\D_\ec$ be the Laplacian associated to $\ec$. Let the potential 
$\fun{\fc}{(-\infty,0]\times S}{\R}$ satisfy hypotheses (F) and (FF) below; let the final condition $\fun{w_0}{S}{\R}$ belong to $\dc(\D_\ec)\cap Lip(S)$ and satisfy (2.13) below (we shall see at the end of section 1 that these conditions are not empty). 

\noindent 1) Then, there is a unique $w\in C^1((-\infty,0],L^2(S,m))\cap C((-\infty,0],\dc(\D_\ec))$ which solves the 
Schr\"odinger equation with time reversed 
$$\left\{
\eqalign{
\partial_t w(t,x)&+\2\D_\ec w(t,x)+\fc(t,x)w(t,x)=0\qquad \forall t\le 0\cr
w(0,x)&=w_0(x)     ,  
}
\right.  \eqno (4)$$
where equalities are in the $L^2(S,m)$ sense, i. e. $m$ a. e..

\noindent 2) The function $w$ is given by the Feynman-Kac formula
$$w(t,x)=E^\ini\left[
e^{\int_t^0\fc(\tau,B_\tau^\ini)\dr\tau}w_0(B_0^\ini)
\right]      $$
for $m$ a. e. $x\in S$. 

\noindent 3) If (2.13) below holds, we shall see that the maximum principle implies 
$w(t,\cdot)$ is bounded away from $0$ and $+\infty$ for all $t\le0$;  we can thus consider  
$$u(t,x)=-\log w(t,x)  .   \eqno (5)$$
Then, $u\in C^1((-\infty,0],L^2(S,m))\cap L^\infty_{loc}((-\infty,0],\dc(\D_\ec))$ and satisfies the Hamilton-Jacobi equation with time reversed 
$$\left\{
\eqalign{
\partial_tu(t,x)&+\2\D_\ec u(t,x)-\2 \G(u(t,\cdot),u(t,\cdot))(x)-\fc(t,x)=0
\qquad \forall t\le 0\cr
u(0,x)&=-\log w_0(x)
}
\right.   \eqno (6)$$
where $\G$ is the carr\'e\ de champs associated with $\ec$ . 

\noindent 4) Lastly, for all probability density $\r_t\in L^\infty(S,m)$ and all $t\le 0$, we have that
$$\int_Su(t,x)\r_t(x)\dr m(x)=
\min\left\{
\int_t^0\dr\tau\int_S\left[
\frac{1}{2}\G(V(\tau,\cdot),V(\tau,\cdot))(y)-\fc(\tau,y)
\right]\dr\mu_\tau(y)+
\int_S u_0(y)\dr\mu_0(y)
\right\}   \eqno (7)$$
where $\mu$ is a solution of the Fokker-Planck equation with drift $V$ starting at 
$\mu_t=\r_t m$; the $\min$ is over all drifts $V\in\wc(\mu)$.

\rm

\vskip 1pc

The paper is organised as follows: in section 1 we recall from [14] and [5] some definitions and theorems about Dirichlet forms; we shall also recall from [2], [3], [4], [5] and [20] the results we need about $RCD(K,\infty)$ spaces. In section 2, we tackle equations (4) and (6) and the Feynman-Kac formula. In section 3, we introduce the notion of weak solutions of the Fokker-Planck equation and prove one inequality of (7); the opposite inequality is proven in section 4.  

\vskip 1pc

\noindent {\bf Acknowledgement.} The author would like to thank the referee for the careful reading and the helpful comments.

\vskip 2pc

\centerline{\bf \S 1}
\centerline{\bf Preliminaries and notation}

\vskip 1pc

To prove points 1) and 2) of theorem 1, it suffices to have a measured metric space with a symmetric, strongly local Dirichlet form on it. We shall call this situation the Dirichlet form setting; let us state the precise hypotheses. 

Following [14], we shall assume that $(S,d)$ is a metric space (which we shall suppose compact for simplicity) and that $m$ is a probability measure on $S$, positive on open sets. 

Let us consider a symmetric bilinear form form $\ec$  
$$\fun{\ec}{\dc(\ec)\times\dc(\ec)}{\R}$$
where $\dc(\ec)$ is a dense subspace of $L^2(S,m)$. The form $\ec$ is called a Dirichlet form if the two conditions below hold. 

\noindent 1)  $\dc(\ec)$ is closed; this means that $\dc(\ec)$ is complete under the norm
$$||u||_{\dc(\ec)}^2=
||u||_{L^2(S,m)}^2+\ec(u,u)  .  $$
It is standard ([14], [17]) that $\ec$ is closed if and only if the quadratic form 
$\fun{}{u}{\ec(u,u)}$ is lower semicontinuous in $L^2(S,m)$. 

\noindent 2)  $\ec$ is Markovian, i. e. 
$$\ec(\eta(f),\eta(f))\le \ec(f,f) $$
for all $f\in\dc(\ec)$ and all 1-Lipschitz maps $\fun{\eta}{\R}{\R}$ with $\eta(0)=0$. 

We shall assume two further properties on $\ec$; the first one is that $\ec$ is regular. This means that $\ec$ has a core, i. e. a subset ${\cal C}\subset\dc(\ec)\cap C(S)$ such that ${\cal C}$ is dense in $\dc(\ec)$ for 
$||\cdot||_{\dc(\ec)}$, and is dense in $C(S)$ for the $\sup$ norm. 

The second one is that $\ec$ is strongly local, i. e. that
$$\ec(f,g)=0$$
if $f,g\in\dc(\ec)$ and $f$ is constant on a neighbourhood of the support of $g$. 

By theorem 1.3.1 of [14], there is a non-positive self-adjoint operator $\D_\ec$ such that 
$$\dc(\D_\ec)\subset\dc(\ec) $$
is dense in $L^2(S,m)$ and
$$\inn{-\D_\ec f}{g}_{L^2(S,m)}=\ec(f,g)\qquad
\forall f\in\dc(\D_\ec),\quad\forall g\in\dc(\ec)  . \eqno (1.1) $$

Now $-\D_\ec$, being self adjoint and non-negative, is monotone maximal; thus we can apply the theory of [10], getting that $-\2\D_\ec$ generates a semigroup of contractions, backward in time, on $L^2(S,m)$. Namely, for $t\le s$ there is 
$\fun{P_{t,s}}{L^2(S,m)}{L^2(S,m)}$ such that
$$P_{t-h-k,t-h}\circ P_{t-h,t}=P_{t-h-k,t}\txt{for} h,k\ge 0  .  $$
The semigroup is autonomous, i. e. $P_{t+h,s+h}=P_{t,s}$, and  we could have called it $P_{t-s}$ as well. The reason for this clumsier notation is that in the next sections we shall need to keep track also of the starting time of the trajectory. 

For each fixed $f\in L^2$ and $t\le s$, the map $\fun{}{t}{P_{t,s}f}$ is continuous and 
$$||P_{t,s}f||_{L^2}\le ||f||_{L^2} .  \eqno (1.2)$$
Since $-\2\D_\ec$ is the generator of $P_{t,s}$, we have that 
$$\2\D_\ec u=\lim_{h\searrow 0}
\frac{P_{0,h}u-u}{h}  =
-\lim_{h\searrow 0}\frac{P_{-h,0}u-u}{-h}  \qquad\forall u\in\dc(\D_\ec)
\eqno (1.3)$$
where the limits are in $L^2(S,m)$.  The Brownian motion is the stochastic process behind the semigroup $P_{t,s}$; namely, by theorem 4.5.3 of [14] for $m$ a. e. 
$x\in S$,  it is possible to define a probability measure ${\bf P}^{(t,x)}$ on $C([t,+\infty),S)$ (and a related expectation $E^{(t,x)}$) such that ${\bf P}^{(t,x)}$ concentrates on 
$$\{
\g\in C([t,+\infty),S)\st\g_t=x
\}   $$
and, for $t\le s$,
$$(P_{t,s}f)(x)=E^{(t,x)}(f(\g_s))   .  \eqno (1.4) $$
For $\tau\in[t,+\infty)$ we define $e_\tau$ as the evaluation map 
$$\fun{e_\tau}{C([t,+\infty),S)}{S},\qquad
\fun{e_\tau}{\g}{\g_\tau}   .   $$
Now the stochastic process 
$$\fun{B^{(t,x)}_\tau}{C([t,+\infty),S)}{S}, \qquad
B^{(t,x)}_\tau\colon=e_\tau\circ\g$$
is the Brownian motion with $B^{(t,x)}_t\equiv x$. 

We denote as usual by $F_\sharp\mu$ the push-forward of a measure $\mu$ by a map 
$F$; for $h>0$ we shall set  
$$p_h(x,\dr y)=(e_h)_\sharp{\bf P}^{(0,x)}  .  $$

By (1.4) we easily get that $P_{t,s}$ is positivity preserving: 
$$P_{t,s}f\ge 0   \txt{if} f\ge 0    .    $$

The last hypothesis we shall make on $\ec$ is that its carr\'e\ de champs is defined on 
$\dc(\ec)$; in other words, we ask that there is a symmetric bilinear form  
$$\fun{\G}{\dc(\ec)\times\dc(\ec)}{L^1(S,m)}$$ 
such that
$$\ec(u,v)=\int_S\G(u,v)\dr m\qquad\forall u,v\in\dc(\ec)  .  \eqno (1.5)$$

We recall from [5] and [9] some properties of the carr\'e\ de champs $\G$. 

\noindent $\bullet)$ If $f,g\in\dc(\ec)$ and $\eta\in Lip(\R)$, then $\eta(f)\in\dc(\ec)$; moreover, the chain rule holds
$$\G(\eta(f),g)=\eta^\prime(f)\G(f,g)   .  \eqno (1.6)$$

\noindent $\bullet$) If $f,g,h\in\dc(\ec)\cap L^\infty$, we have the Leibnitz rule
$$\G(fg,h)=f\G(g,h)+g\G(f,h)  .  \eqno (1.7)$$

\noindent $\bullet$) By (1.5), (1.1) becomes
$$\inn{-\D_\ec u}{v}_{L^2(S,m)}=
\int_S \G(u,v)\dr m\qquad 
\forall u\in\dc(\D_\ec),\quad\forall v\in\dc(\ec)  \eqno (1.8)$$

\noindent $\bullet$) By formula (2.18) of [5] we have that, if $\eta\in C^2(\R)$ has bounded first and second derivatives (the hypothesis $\eta(0)=0$ is not necessary when $m$ is finite), if $f\in\dc(\D_\ec)$ and $\G(f,f)\in L^2(S,m)$, then 
$\eta(f)\in L^2(S,m)$ and
$$\D_\ec\eta(f)=\eta^\prime(f)\D_\ec f+
\eta^\pprime(f)\G(f,f)  .  \eqno (1.9)$$

There is an important example of a Dirichlet form which satisfies all these hypotheses: the double of Cheeger's energy on $RCD(K,\infty)$ spaces. When we shall be in this more restrictive hypothesis, we shall say that we are in the $RCD(K,\infty)$ setting. We refer the reader to [2], [3], [4], [5] and [20] for their definition and the study of their properties; here, we only recall the few facts we need. 

\noindent $\bullet$) There is a "natural" Dirichlet form $\ec$;  $\ec$ is the double of Cheeger's energy, which we don't define (see for instance [2].) The form $\ec$ is regular, strongly local and its carr\'e\ de champs is defined on $\dc(\ec)$. 

\noindent $\bullet$) Each $x\in S$ is the starting point of a Brownian motion.

\noindent $\bullet$) We recall that we called $p_h(x,\dr y)$ the transition probability of the Brownian motion. We define $\pc(S)$ as the space of all the Borel probability measures on $S$; if $\mu,\nu\in\pc(S)$, we define their Wasserstein distance $W_2(\mu,\nu)$ as
$$W_2^2(\mu,\nu)=\min\int_{S\times S}d^2(x,y)\dr\Sigma(x,y)   $$
where the minimum is over all Borel probability measures $\Sigma$ on $S\times S$ whose first and secon marginals are, respectively, $\mu$ and $\nu$. 

An important property of $RCD(K,\infty)$ spaces is that the map $\fun{}{x}{p_h(x,\dr y)}$ is Lipschitz from $(S,d)$ to $(\pc(S),W_2)$; the Lipschitz constant is bounded by $e^{-Kh}$. 

\noindent $\bullet$) The probability $p_h(x,\dr y)$ has a density: 
$p_h(x,\dr y)=\check p_h(x,y)m(\dr y)$. The function $\check p_h$ is symmetric, i. e. 
$\check p_h(x,y)=\check p_h(y,x)$

\noindent $\bullet$) If $\fun{a}{S}{\R}$ is Lipschitz, then
$$\G(a,a)\le (Lip(a))^2   .   \eqno (1.10)$$

\noindent $\bullet$) By section 4.1 of [5], Lipschitz functions are dense in $\dc(\D_\ec)$ for the $L^2$ topology. As a consequence, in theorem 1 the hypotheses on $w_0$ are not empty.

\vskip 2pc
\centerline{\bf \S 2}
\centerline{\bf The Feynman-Kac formula}
\vskip 1pc

In this section, we are going to prove that the Schr\"odinger and Hamilton-Jacobi equations have a unique solution and that the Feynman-Kac formula holds; as usual, Feynman-Kac will imply a maximum principle. We shall also use the Feynman-Kac formula to prove that the solution of Schr\"odinger's equation (4) is Lipschitz; this will help us to deduce equation (6) at the end of this section. 

We start in the Dirichlet form setting. We saw above that $-\2\D_\ec$ is a monotone maximal operator; in particular, its graph is closed. Said differently, $\dc(\D_\ec)$ with the internal product 
$$\inn{u}{v}_{\D_\ec}\colon=
\inn{u}{v}_{L^2(S,m)}+
\inn{\D_\ec u}{\D_\ec v}_{L^2(S,m)}    $$
is a Hilbert space. We shall denote its norm by 
$$||u||^2_{\D_\ec}\colon=||u||^2_{L^2}+||\D_\ec u||^2_{L^2}  .  $$

Let $g\in C((-\infty,0],L^2(S,m))$; shall say that $u$ is a strong solution of the inhomogeneous heat equation, backward in time, 
$$\left\{
\eqalign{
\partial_tu&=-\2\D_\ec u-g(t,\cdot)\quad t\le 0\cr
u(0)&=u_0
}
\right.   \eqno (2.1)$$
if $u\in C^1([-\infty,0],L^2(S,m))\cap C([-\infty,0],\dc(\D_\ec))$ and if (2.1) holds for all 
$t\in(-\infty,0]$; the equalities are in $L^2$, i. e. m a. e.. According to [16], theorem IX.1.19 (or corollary 4.2.5 of [18]), (2.1) has a unique strong solution if $u_0\in\dc(\D_\ec)$ and 
$g\in C^1((-\infty,0], L^2)$. Moreover, the solution is given by the Duhamel formula 
$$u(t,\cdot)=P_{t,0}u_0+
\int_t^0P_{t,s}g(s,\cdot)\dr s, \qquad t\le 0  .   \eqno (2.2)$$
Let us set 
$$||u||_{C^1([-T,0],L^2)}=\sup_{t\in[-T,0]}
[
||u(t,\cdot)||_{L^2}+ 
||\partial_t u(t,\cdot)||_{L^2} 
]    .   $$

From now on we shall make the following assumption on the potential of the Schr\"odinger equation $\fc$:

\vskip 1pc

\noindent{\bf(F)}:  $\fc\in C^1((-\infty,0],C(S,\R))$. 

\vskip 1pc

The following lemma is well-known (see for instance theorems 6.1.4 and 6.1.5 of [18]); we give the proof for completeness. 

\lem{2.1} Let $(S,d,m)$ and $\ec$ satisfy the hypotheses of the Dirichlet form setting. Let 
$\fc$ satisfy assumption (F) and let $w_0\in\dc(\D_\ec)$. Then, the 
Schr\"odinger equation with time reversed
$$\left\{
\eqalign{
\partial_t w(t,x)&+\2\D_\ec w(t,x)+\fc(t,x)w(t,x)=0\quad t\le 0\cr
w(0,x)&=w_0(x)
}
\right.   \eqno (2.3)$$
has a unique strong solution
$$w\in C^1((-\infty,0],L^2)\cap C((-\infty,0],\dc(\D_\ec))   .   \eqno (2.4)$$

\proof By hypothesis (F), if $w\in C^1((-\infty,0],L^2)$, then also 
$\fc w\in C^1((-\infty,0],L^2)$; by the theorem of [16] we mentioned before formula (2.2), we see that it suffices to find a solution $w\in C^1((-\infty,0],L^2)$ of the integral equation
$$w(t,\cdot)=P_{t,0}w_0+
\int_t^0P_{t,s}[\fc(s,\cdot)w(s,\cdot)]\dr s   .     $$
Setting $s=t-r$, this becomes (note that $r\le 0$) 
$$w(t,\cdot)=
P_{t,0}w_0+
\int_t^0 P_{r,0}[
\fc(t-r)w(t-r)
]   \dr r    .  \eqno (2.5)$$ 
For $T>0$ let us set 
$$A_T=\{
w\in C^1([-T,0],L^2)\st w(0,\cdot)=w_0
\}   .   $$
Clearly, $A_T$ is closed in $C^1([-T,0],L^2)$. We begin to show that, for $T>0$ small, the operator
$$\fun{\Phi}{A_T}{A_T}$$
$$\fun{\Phi}{w}{
P_{t,0}w_0+
\int_t^0 P_{t,s}[\fc(s,\cdot)w(s,\cdot)] \dr s
}   $$
is a contraction from $A_T$ to itself and thus it has a unique fixed point. 

Formula (2.5) and an easy calculation imply that 
$$\partial_t\Phi(w)(t,\cdot)=-\2\D_\ec P_{t,0}w_0
-P_{t,0}[
\fc(0,\cdot)w(0,\cdot)
]  +
\int_t^0 P_{r,0}\left\{
\frac{{\rm d}}{{\rm d}t}[
\fc(t-r,\cdot)w(t-r,\cdot)
]
\right\}   \dr r   .  $$
If $w,\tilde w\in A_T$, this yields the equality below; the first inequality comes from H\"older and (1.2); the last inequality follows since $t\in[-T,0]$. 
$$||\partial_t\Phi(w)(t)-\partial_t\Phi(\tilde w)(t)||_{L^2}=$$
$$\left\vert\left\vert
\int_t^0 P_{r,0}\left\{
\frac{\rm d}{{\rm d}t}[
\fc(t-r)w(t-r)-\fc(t-r)\tilde w(t-r)
]
\right\}  \dr r
\right\vert\right\vert_{L^2}\le$$
$$\int_t^0||\fc||_{C^1([-T,0],C(S))}\cdot
||w-\tilde w||_{C^1([-T,0],L^2)}\dr r\le$$
$$T||\fc||_{C^1([-T,0],C(S))}\cdot
||w-\tilde w||_{C^1([-T,0],L^2)}   .   $$
Thus, if $T$ is so small that
$$T||\fc||_{C^1([-T,0],L^2)}\le\frac{1}{4} ,  \eqno (2.6)$$
we get that
$$||\partial_t[\Phi(w)-\Phi(\tilde w)]||_{C^0([-T,0],L^2)}\le
\frac{1}{4}||w-\tilde w||_{C^1([-T,0],L^2)} . $$
Analogously, 
$$||[\Phi(w)-\Phi(\tilde w)]||_{C^0([-T,0],L^2)}\le
\frac{1}{4}||w-\tilde w||_{C^1([-T,0],L^2)} . $$
By the last two formulas, we see that $\Phi$ is a contraction. 

Deducing from this existence for all times is standard. It suffices to show that there is a decreasing function $\fun{\e}{[0,+\infty)}{[0,+\infty)}$ independent of $u$ and such that, if the solution exists up to time $-T$, then it can be extended to $-T-\e(T)$. This implies in a standard way that the maximal interval of existence of $w$ is $(-\infty,0]$. 

We find the function $\e$. We choose $w(-T)$ as a final condition for the operator $\Phi$; in other words, we set
$$\tilde A_{[-T-\e,-T]}=\{
\hat w\in C^1([-T-\e,-T],L^2)\st\hat w(-T,\cdot)=w(-T,\cdot)
\}    $$
and define 
$$\fun{\tilde\Phi}{A_{[-T-\e,-T]}}{A_{[-T-\e,-T]}}  $$
$$\fun{\tilde\Phi}{\hat w}{P_{t,-T}\hat w(-T,\cdot)+\int_t^{-T}P_{t,s}[\fc(s,\cdot)
\hat w(s,\cdot)]\dr s}  .  $$
Arguing as above, we see that the Lipschitz constant of $\tilde\Phi$ is smaller than 
$\frac{1}{2}$ as long as   
$$\e ||\fc||_{C^1([-T-\e,-T],L^\infty)}\le\frac{1}{4}  .  $$
Now it suffices to take $\e(T)$ as the largest $\e$ for which the formula above holds. 

\fin

Next, we prove the Feynman-Kac formula.

\lem{2.2} Let $(S,d,m)$ and $\ec$ satisfy the hypotheses of the Dirichlet form setting. Let 
$\fc$ satisfy (F), let $w_0\in\dc(\D_\ec)$ and let $w$ be the unique solution of (2.3). Then, for $m$ a. e. $x\in S$ we have that
$$w(t,x)=E^\ini\left[
e^{
\int_t^0\fc(\tau,B^\ini_\tau)\dr\tau
}   w(0,B_0^\ini)
\right]    .   \eqno (2.7)$$

\proof We recall the argument of [11]. We fix $t<0$; for 
$s\in[t,0]$ we set 
$$a(s,x)=E^\ini\left[
e^{
\int_t^s\fc(\tau,B^\ini_\tau)\dr\tau
} w(s,B_s^\ini)
\right]   .   \eqno (2.8)$$
We are going to prove that, for all $\psi\in L^2$, the function
$$\fun{l_\psi}{s}{\inn{a(s,\cdot)}{\psi}_{L^2}}  $$
has the right hand derivative identically equal to zero at every $s\in[t,0)$. It is a standard fact (see for instance [19], theorem 8.21) that this implies that $l_\psi$ is absolutely continuous on $(t,0)$. Using (2.8) and dominated convergence, we easily see that $l_\psi$ is continuous at $s=0$. Thus, integrating $l_\psi^\prime$, we get that
$$\inn{a(0,\cdot)}{\psi}_{L^2}=\inn{a(t,\cdot)}{\psi}_{L^2}
\qquad\forall\psi\in L^2  .  $$
Thus, $a(0,\cdot)=a(t,\cdot)$ in $L^2$, which implies (2.7). 

We calculate the right hand derivative for $l_\psi$. The first equality below comes from the definition of $a(s,x)$ in (2.8). 
$$\frac{\dr^+}{\dr s}\inn{a(s,\cdot)}{\psi}_{L^2}=$$
$$\lim_{h\searrow 0}\inn{
\frac{1}{h}E^\ini\left[
e^{
\int_t^{s+h}\fc(\tau,B^\ini_\tau)\dr\tau
}  w(s+h,B^\ini_{s+h})-
e^{
\int_t^{s}\fc(\tau,B^\ini_\tau)\dr\tau
}  w(s,B^\ini_{s})
\right]
}{\psi}_{L^2}  =  $$
$$\lim_{h\searrow 0}\inn{
\frac{1}{h}E^\ini\left\{
\left[
e^{
\int_t^{s+h}\fc(\tau,B^\ini_\tau)\dr\tau
} -e^{
\int_t^{s}\fc(\tau,B^\ini_\tau)\dr\tau
}
\right]  w(s+h,B^\ini_{s+h})
\right\}
}{\psi}_{L^2}+    \eqno(2.9)_a$$
$$\lim_{h\searrow 0}\inn{
\frac{1}{h}E^\ini\left\{
e^{
\int_t^{s}\fc(\tau,B^\ini_\tau)\dr\tau
}    \left[
w(s+h,B^\ini_{s+h})-w(s,B^\ini_{s})
\right]
\right\}
}{\psi}_{L^2}  .  \eqno (2.9)_b$$
We begin to tackle $(2.9)_a$. 
We consider the measure $m\otimes{\bf P}^{(t,x)}$ on $S\times C([t,0],S)$; since $\fc$ is continuous and the Brownian motion has continuous trajectories, we have that 
$$\frac{1}{h}
\left[
e^{
\int_t^{s+h}\fc(\tau,B^\ini_\tau)\dr\tau
} -e^{
\int_t^{s}\fc(\tau,B^\ini_\tau)\dr\tau
}
\right]  \tends
e^{
\int_t^{s}\fc(\tau,B^\ini_\tau)\dr\tau
}
\fc(s,B^{(t,x)}_s)$$
$m\otimes{\bf P}^{(t,x)}$ a. e.; adding in the fact that $\fc$ is bounded, we get that 
$$\left\vert
\frac{1}{h}
\left[
e^{
\int_t^{s+h}\fc(\tau,B^\ini_\tau)\dr\tau
} -e^{
\int_t^{s}\fc(\tau,B^\ini_\tau)\dr\tau
}
\right]
\right\vert     \le M  \eqno (2.10)$$
$m\otimes{\bf P}^{(t,x)}$ a. e.. By dominated convergence, this implies that, for all 
$h\in(0,1]$, 
$$E^{(t,x)}\left\{
\frac{1}{h}
\left[
e^{
\int_t^{s+h}\fc(\tau,B^\ini_\tau)\dr\tau
} -e^{
\int_t^{s}\fc(\tau,B^\ini_\tau)\dr\tau
}
\right] w(s,B_s^{(t,x)})
\right\}   \tends$$
$$E^{(t,x)}\left[
e^{
\int_t^{s}\fc(\tau,B^\ini_\tau)\dr\tau
}
\fc(s,B^{(t,x)}_s) w(s,B^{(t,x)}_s)
\right]     \eqno (2.11)$$
in $L^2(S,m)$.

We assert that
$$E^{(t,x)}[
w(s+h,B_{s+h}^{(t,x)})-w(s,B^{(t,x)}_s)
]\tends 0   \txt{in}  L^2(S,m)  .  \eqno (2.12)$$
We prove this; the first equality below is the definition of $P_{t,s+h}$; the inequality is (1.2) and the limit follows from the regularity of $w$.
$$||
E^{(t,x)}[w(s+h,B_{s+h}^{(t,x)})-w(s,B_{s+h}^{(t,x)})]
||_{L^2}=
||
P_{t,s+h}[w(s+h,\cdot)-w(s,\cdot)]
||_{L^2}\le$$
$$||
w(s+h,\cdot)-w(s,\cdot)
||_{L^2}  \tends 0   \txt{as}h\tends 0.  $$
The first equality below is the definition of $P_{t,s}$, while the limit comes from the fact that $P_{t,s}$ is strongly continuous. 
$$||
E^{(t,x)}[w(s,B_{s+h}^{(t,x)})-w(s,B_{s}^{(t,x)})]
||_{L^2}   =
||
[P_{t,s+h}-P_{t,s}]w(s,\cdot)
||_{L^2}  \tends 0   \txt{ as $h\tends 0$.}  $$
Now (2.12) follows from the last two inequalities. 
 
By (2.10) and (2.12) we easily get that
$$\frac{1}{h}E^\ini\left\{
\left[
e^{
\int_t^{s+h}\fc(\tau,B^\ini_\tau)\dr\tau
} -e^{
\int_t^{s}\fc(\tau,B^\ini_\tau)\dr\tau
}
\right] [ w(s+h,B^\ini_{s+h})-w(s,B^{(t,x)}_s)]
\right\}
\tends   0     $$
in $L^2(S,m)$. 

Together with (2.11), this implies the convergence below, which is in the $L^2$ topology. 
$$\frac{1}{h}E^\ini\left\{
\left[
e^{
\int_t^{s+h}\fc(\tau,B^\ini_\tau)\dr\tau
} -e^{
\int_t^{s}\fc(\tau,B^\ini_\tau)\dr\tau
}
\right]  w(s+h,B^\ini_{s+h})
\right\}   =$$
$$\frac{1}{h}E^{(t,x)}\left\{
\left[
e^{\int_t^{s+h}\fc(\tau,B^{(t,x)}_\tau)\dr\tau}-
e^{\int_t^{s}\fc(\tau,B^{(t,x)}_\tau)\dr\tau}
\right]
[
w(s+h,B^{(t,x)}_{s+h})-w(s,B^{(t,x)}_{s})
]
\right\}   +$$
$$\frac{1}{h}E^{(t,x)}\left[
e^{
\int_t^{s}\fc(\tau,B^\ini_\tau)\dr\tau
}
\fc(s,B^{(t,x)}_s) w(s,B^{(t,x)}_s)
\right]
\tends$$
$$E^\ini\left[
e^{
\int_t^{s}\fc(\tau,B^\ini_\tau)\dr\tau
}\fc(s,B^\ini_s)w(s,B^\ini_s)
\right]   . $$
 
The formula above yields the limit in $(2.9)_a$; for the limit in $(2.9)_b$ we argue analogously. Using  (1.3), the fact that 
$w\in C^1((-\infty,0],L^2)\cap C^0((-\infty,0],\dc(\D_\ec))$ and the Markovianity of the Brownian motion imply that 
$$\lim_{h\searrow 0}
\frac{1}{h}E^\ini\left\{
e^{
\int_t^{s}\fc(\tau,B^\ini_\tau)\dr\tau
}    \left[
w(s+h,B^\ini_{s+h})-w(s,B^\ini_{s})
\right]
\right\}
=$$
$$
E^\ini\left\{
e^{
\int_t^{s}\fc(\tau,B^\ini_\tau)\dr\tau
}
[
\partial_s w(s,B^\ini_s)+\2\D_\ec w(s,B^\ini_s)
]
\right\}   $$
where the limit is in $L^2(S,m)$. 

The first equality below comes from (2.9) and the last two formulas, the second one comes from the fact that $w$ solves (2.3). 
$$\frac{\dr^+}{\dr s}\inn{a(s,\cdot)}{\psi}_{L^2}=$$
$$\inn{
E^{(t,x)}\left\{
e^{
\int_t^s\fc(\tau, B^{(t,x)}_\tau)\dr\tau
}
\left[
\partial_s w(s,B^{(t,x)}_s)+\2\D_\ec w(s,B^{(t,x)}_s)+\fc(s,B^{(t,x)}_s)w(s,B^{(t,x)}_s)
\right]
\right\}
}{\psi}   =0   $$
$$\forall s\in(t,0)$$
as we wanted.

\fin

Note that the integral in (2.7) converges also if the initial condition $w(0,\cdot)$ is only bounded; actually, we shall show in lemma 2.4 below that, if $w_0$ is Lipschitz, also 
$w(t,\cdot)$ is such. Naturally, if $w_0\not\in\dc(\D_\ec)$, we lose the fact that $w$ solves (2.3). 

An immediate consequence of (2.7) and hypothesis (F) is the maximum principle. 

\cor{2.3} Let $(S,d,m)$ and $\ec$ satisfy the hypotheses of the Dirichlet form setting. Let 
$\fc$ satisfy (F), let $w_0\in C(S)$, let $w$ be defined as in (2.7) and let us suppose that, for some $c_1>0$, 
$$\frac{1}{c_1}\le w_0\le c_1    .    \eqno (2.13)$$
Then, there is an increasing function $\fun{D_1}{[0,+\infty)}{[0,+\infty)}$ such that
$$\frac{1}{D_1(T)}\le w(t,\cdot)\le D_1(T)
\txt{for all}t\in[-T,0]  \txt{and $m$ a. e.} x\in S  .  \eqno (2.14) $$

\rm

\vskip 1pc

Up to this point, we have only used the properties of a strongly local Dirichlet form in a metric space; if we want to prove formula (6) of the introduction, we need much more. Thus, from now on we strengthen our hypotheses to the $RCD(K,\infty)$ setting. We shall need a stronger condition on $\fc$ too. 

\vskip 1pc

\noindent {\bf(FF)} $\fc(t,\cdot)$ is Lipschitz for all $t\le 0$ and there is an increasing function $\fun{D_2}{[0,+\infty)}{[0,+\infty)}$ such that
$$Lip(\fc(t,\cdot))\le D_2(T)
\txt{if} t\in[-T,0]. $$

\vskip 1pc

\lem{2.4} Let $(S,d,m)$ and $\ec$ satisfy the hypotheses of the $RCD(K,\infty)$ setting.  Let $w$ be given by (2.7), with $w(0,\cdot)$ Lipschitz; let $\fc$ satisfy (F) and 
(FF). Then, there is an increasing function $\fun{D_3}{[0,+\infty)}{[0,+\infty)}$ such that $w(t,\cdot)$ is 
$D_3(T)$-Lipschitz for all $t\in[-T,0]$. 

\proof Since $w$ is defined by (2.7), hypothesis (F) and dominated convergence give the first equality below; the second one is the definition of Wiener's measure. 
$$w(t,x)=\lim_{n\tends+\infty} E^{(t,x)}\left[
\exp\left(
\frac{1}{n}\sum_{j=1}^n\fc(t(1-\frac{j}{n}),B^{(t,x)}_{t(1-\frac{j}{n})})
\right)
w_0(B_0^{(t,x)})
\right]   =   $$
$$\lim_{n\tends+\infty}
\int_S\exp\left(
\frac{1}{n}\fc(t\left(1-\frac{1}{n}\right),x_{t\left(1-\frac{1}{n}\right)})
\right)
p_\frac{|t|}{n}(x,\dr x_{t(1-\frac{1}{n})})$$
$$\int_S\exp\left(
\frac{1}{n}\fc(t\left(1-\frac{2}{n}\right),x_{t\left(2-\frac{1}{n}\right)})
\right)
p_\frac{|t|}{n}(x_{t(1-\frac{1}{n})},\dr x_{t(1-\frac{2}{n})})  \dots$$
$$\dots\int_S
\exp\left(
\frac{1}{n}\fc\left(0,x_0\right)
\right)
w_0(x_0)
p_\frac{|t|}{n}(x_\frac{t}{n},\dr x_{0})   .    $$
Thus, it suffices to prove that the multiple integral on the right is $D_3(T)$-Lipschitz in $x$ for all $t\in[-T,0]$ and all $n\ge 1$. By the particular form of the integral above, this follows 
by iteration if we show that, when $a$ is Lipschitz, the map 
$$\fun{b}{x}{
\int_S\exp\left(\frac{1}{n}\fc(\frac{j}{n},z)\right)  a(z)p_\frac{|t|}{n}(x,\dr z)
}  $$
satisfies
$$Lip(b)\le e^\frac{D_4(T)}{n}[Lip(a)+||a||_\infty]
\txt{and}
||b||_\infty\le e^\frac{D_4(T)}{n}||a||_\infty  .  $$
The inequality on the right comes easily from the definition of $b$ and hypothesis (F); let us prove the one on the left. We recall from section 1 that the map 
$\fun{}{x}{p_\frac{|t|}{n}}(x,\dr z)$ from $S$ to $\pc(S)$ is $e^\frac{-K}{n}$-Lipschitz; together with the fact that
$$Lip\left[
\exp\left(
\frac{1}{n}\fc\left(
\frac{j}{n},z
\right)
\right)   a(z)
\right]   \le
e^\frac{D_5(T)}{n}[Lip(a)+||a||_\infty]    $$
(which follows from (FF) and the formula for the Lipschitz constant of a product) this implies the inequality on the left. 

\fin

\lem{2.5} Let $(S,d,m)$ and $\ec$ satisfy the hypotheses of the $RCD(K,\infty)$ setting.   Let $\fc$ satisfy (F) and (FF), let $w_0\in\dc(\D_\ec)\cap Lip(S)$ satisfy (2.13); let $w$ be the solution of (2.3). By corollary 2.3 we can set 
$$u(t,x)=-\log w(t,x)  .  $$
Then, $u\in C^1((-\infty,0], L^2(S,m))\cap L^\infty_{loc}((-\infty,0],\dc(\D_\ec))$ and solves the Hamilton-Jacobi equation 
$$\left\{
\eqalign{
\partial_t u(t,\cdot) &+\2\D_\ec u(t,\cdot)-\2 \G(u(t,\cdot),u(t,\cdot))
-\fc(t,\cdot)=0   \quad t\le 0\cr
u(0,\cdot)&=u_0\colon=-\log w_0    
}
\right.$$
where the equalities are in $L^2(S,m)$, i. e. $m$ a. e..

\proof Lemma 2.1 and corollary 2.3 imply that $u\in C^1((-\infty,0],L^2)$. By corollary 2.3, 
$\eta(x)=-\log(x)$ has bounded first and second derivatives on the range of $w$, and $\eta\circ w\in L^2(S,m)$. By lemma 2.4 and (1.10), we can apply (1.9) with $f=w$ and we get the first equality below, while the second one comes from (2.3); the last one comes from the chain rule (1.6). 
$$\2\D_\ec u(t,x)=
-\frac{1}{w(t,x)}\cdot\2\cdot\D_\ec w(t,x)+
\2\cdot\frac{1}{w^2(t,x)}\cdot \G(w(t,\cdot), w(t,\cdot))(x)=$$
$$\frac{1}{w(t,x)}  [
\partial_t w(t,x)+\fc(t,x)\cdot w(t,x)
]    
+\2\cdot\frac{1}{w^2(t,x)}\cdot \G(w(t,\cdot), w(t,\cdot))(x)=$$
$$-\partial_tu(t,x)+\2\G(u(t,\cdot), u(t,\cdot))(x)+\fc(t,x)  .  $$

\fin

\vskip 2pc
\centerline{\bf \S 3}
\centerline{\bf The weak version of Fokker-Planck and the value function}
\vskip 1pc

In this section, we define the weak version of the Fokker-Planck equation and prove one inequality of (7). In this section and in the next one, we shall suppose that $(S,d,m)$ is a 
$RCD(K,\infty)$ space and that $\ec$ is the double of Cheeger's energy. 

\noindent{\bf The space of test functions.} We recall from section 2 that $\dc(\D_\ec)$ with the internal product $\inn{\cdot}{\cdot}_{\D_\ec}$ is a Hilbert space. 

Let $u\in\dc(\D_\ec)$; since $\dc(\D_\ec)\subset\dc(\ec)$, we can take $g=u$ in (1.8) getting that 
$$\int_S\G(u,u)\dr m\le
||\D_\ec u||_{L^2}\cdot ||u||_{L^2}  .  $$
By the definition of $||\cdot||_{\dc(\ec)}$ and Young's inequality, this implies that  
$$||u||_{\dc(\ec)}\le
\sqrt\frac{3}{2}\cdot ||u||_{\D_\ec}\qquad
\forall u\in\dc(\D_\ec)  .  $$
In particular, we have that, if $t<0$, 
$$L^\infty([t,0],\dc(\ec))\subset
L^\infty([t,0],\dc(\D_\ec))  .  \eqno (3.1)$$

We say that $\phi$ is a test function, or that $\phi\in\tc$ for short, if
$$\phi\in C^1([t,0],L^2)\cap L^\infty([t,0],\dc(\D_\ec))\cap
L^\infty([t,0]\times S,\L^1\otimes m)  .  $$
By (3.1), we have that
$$\tc\subset  L^\infty([t,0],\dc(\ec))  .    \eqno (3.2)$$

\vskip 1pc
\noindent{\bf The space of drifts.} If $\fun{\mu}{[t,0]}{\pc(S)}$ is a curve of measures, we shall need to compare the "tangent spaces" to $\pc(S)$ at $\mu_\tau$ and 
$\mu_{\tau^\prime}$; as we shall see below, the following definition radically simplifies the problem. 

We say that a Borel function $\fun{\mu}{[t,0]}{\pc(S)}$ is admissible if $\mu_\tau=\r_\tau m$ for all $\tau\in[t,0]$ and if there is $C_1>0$ such that
$$\r_\tau\le C_1  \qquad\forall \tau\in[t,0]    .    \eqno (3.3)$$
Let $\fun{\mu}{[t,0]}{\pc(S)}$ be admissible; we define $\wc(\mu)$ as the space of the Borel functions $\fun{u}{(t,0)}{\dc(\ec)}$ such that 
$$||u||^2_{\wc(\mu)}\colon=
\int_t^0\dr \tau\int_S|u(\tau,x)|^2\dr\mu_\tau(x)+
\int_t^0\dr \tau\int_S\G(u(\tau,\cdot),u(\tau,\cdot))(x)\dr\mu_\tau(x)<+\infty  .  $$
Naturally, an important question is whether the test functions $\tc$ are dense in $\wc(\mu)$ (actually, in [1] 
$\wc(\mu)$ is {\it defined} as the closure of the gradients of the test functions in the suitable topology); we don't address this question because we are not going to need the answer. We shall only need the trivial fact that $\tc\subset\wc(\mu)$. 

\noindent{\bf The equation.} Let $\fun{\mu}{[t,0]}{\pc(S)}$ be Borel and satisfy (3.3); we say that $\mu$ is a weak solution of the Fokker-Planck equation with drift $V\in\wc(\mu)$ if
$$\int_t^0\dr \tau\int_S[
\partial_\tau\phi_\tau+\2\D_\ec\phi_\tau+\G(V_\tau,\phi_\tau)
]  \dr\mu_\tau=$$
$$\int_S\phi_0(x)\dr\mu_0(x)-\int_S\phi_t(x)\dr\mu_t(x)   \qquad
\forall\phi\in\tc  . \eqno (3.4)$$

We check that the integral on the left makes sense. Note that, since $\phi\in\tc$, the second inequality below follows, while the first one is H\"older. 
$$\sup_{\tau\in[t,0]}||\partial_\tau\phi||_{L^1(m)}\le 
\sup_{\tau\in[t,0]}||\partial_\tau\phi||_{L^2(m)}<+\infty . $$
Now (3.3) implies that $\partial_\tau\phi\in L^1([t,0]\times S,\L^1\otimes\mu_\tau)$. To prove that  
$\G(V_\tau,\phi_\tau)\in L^1([t,0]\times S,\L^1\otimes\mu_\tau)$ we note that, by 
Cauchy-Schwarz,
$$|\G(\phi_t,V_t)|\le
\G(\phi_t,\phi_t)^\2\cdot\G(V_t,V_t)^\2  \le
\2[
\G(\phi_t,\phi_t)+\G(V_t,V_t)
]  .   $$
Now $\G(V_\tau,V_\tau)\in L^1(\L^1\otimes\mu_\tau)$ because 
$V\in\wc(\mu)$, and $\G(\phi_\tau,\phi_\tau)\in L^1(\L^1\otimes\mu_\tau)$ since 
$\phi\in\tc$ and (3.2), (3.3) hold; thus, the last inequality implies that 
$\G(V_\tau,\phi_\tau)\in L^1([t,0]\times S,\L^1\otimes\mu_\tau)$.  

Since $\phi\in\tc$, the map $\fun{}{\tau}{\D_\ec\phi_\tau}$ is bounded from 
$[t,0]$ to $L^2(S,m)$; together with (3.3), this implies the integrability of 
$\2\D_\ec\phi_\tau$.

As a last remark, note that we don't address the question whether a Borel curve of measures satisfying (3.4) is continuous or absolutely continuous; we only note that the solution of Fokker Planck we build in section 4 is continuous. We refer the reader to [12] for a study of this problem on $\R^d$. 

\vskip 1pc

\noindent{\bf Beginning of the proof of theorem 1.} The solution of (4) exists and is unique by lemma 2.1; the solution of (6), by lemma 2.5. The Feynman-Kac formula of point 2) follows from lemma 2.2. Thus, we are left with proving (7); the following lemma gives one of the inequalities.

\lem{3.1} Let $(S,d,m)$ satisfy the hypotheses of the $RCD(K,\infty)$ setting and let (F) and (FF) hold. Let $u$ be defined as in lemma 2.5; let $t<0$ and let $\mu_\tau=\r_\tau m$ be a Borel curve of measures which satisfies (3.3) and the Fokker-Planck equation (3.4) on $[t,0]$ for a drift 
$V\in\wc(\mu)$. Then, 
$$\int_S u(t,x)\r_t(x)\dr m(x)\le
\int_t^0\dr\tau
\int_S[
\2\G(V(\tau,\cdot),V(\tau,\cdot))(x)-\fc(\tau,x)
]\dr\mu_\tau(x)+
\int_Su_0\dr\mu_0  .  \eqno (3.5)$$
Equality holds in the formula above when $V=-u$. 

\proof We note that $u\in\tc$: indeed, it belongs to 
$C^1([t,0],L^2)\cap L^\infty([t,0],\dc(\D_\ec))$ by lemma 2.5, and to $L^\infty([t,0]\times S)$ by corollary 2.3. Thus, we can use $u$ as a test function in the Fokker-Planck equation (3.4) on the time-interval $[t,0]$ and get the first equality below. The second one comes from the fact that $u$ solves (6) and the fact that, by (3.3), $m$ a. e. implies $\mu_t$ a. e.. The last inequality comes from the properties of quadratic forms. 
$$\int_S u(t,x)\r_t(x)\dr m(x)=
\int_S u(0,y)\dr\mu_{0}(y)+$$
$$\int_t^{0}\dr\tau\int_S[
-\partial_\tau u(\tau,y)-\2\D_\ec u(\tau,y)-\G(V,u(\tau,\cdot))(y)
]  \dr\mu_\tau(y) = $$
$$\int_S u(0,y)\dr\mu_{0}(y)+
\int_t^{0}\dr\tau\int_S[
-\2 \G(u(\tau,\cdot),u(\tau,\cdot))(y)-\G(V,u(\tau,\cdot))(y)-\fc(\tau,y)
]  \dr\mu_\tau(y)\le$$
$$\int_S  u_0(y)\dr\mu_{0}(y)+
\int_t^{0}\dr\tau\int_S[
\2\G(V(\tau,\cdot),V(\tau,\cdot))(y)-\fc(\tau,y)
]  \dr\mu_\tau(y)    .  $$

If $V=-u$, then the only inequality in the formula above becomes an equality, implying the last assertion of the lemma. 

\fin

For the opposite inequality, we need to solve the Fokker-Planck equation; that's what we do in the next section. 

\rm

\vskip 2pc
\centerline{\bf \S 4}
\centerline{\bf Solving Fokker-Planck}
\vskip 1pc

From the last assertion of lemma 3.1 we deduce that the following proposition implies the inequality opposite to (3.5).  

\prop{4.1} Let $(S,d,m)$ satisfy the hypotheses of the $RCD(K,\infty)$ setting and let (F) and (FF) hold. Let $u$ be as in lemma 2.5; let $t<0$ and let $\r_t$ be a bounded probability density on $S$. Then, there is a continuous curve $\fun{\mu}{[t,0]}{\pc(S)}$ such that

\noindent 1) $\mu$ is admissible, i. e. it satisfies (3.3) and

\noindent 2) $\mu$ is a weak solution of the Fokker-Planck equation with drift $-u$ and initial condition 
$\mu_t=\r_t m$. 

\rm

\vskip 1pc

The idea of the proof is the following: if $\fun{}{s}{\mu^{(t,x)}_s}$ is the flow of the Fokker-Planck equation with drift $-u$, starting at $\mu^{(t,x)}_t=\d_x$, $t\le s$, and $\b\in C(S)$, then
$$f(t,x)\colon =\int_S \b(y)\dr\mu^{(t,x)}_s(y)$$
is (at least morally) a solution of the conjugate of Fokker-Planck on $(-\infty,s)$ with 
$f(s,\cdot)=\b$. We are going to work our way backwards: first we show that the conjugate of Fokker-Planck has a solution; then, we  define $\mu^{(t,x)}_s$ by the formula above and the Riesz representation theorem. Point 1) of proposition 4.1 will follow by some estimates on $f$; point 2) will follow from the fact that $f$ satisfies the conjugate of Fokker-Planck.

We shall need the stronger hypothesis that $(S,d,m)$  is a $RCD(K,\infty)$ space because we want to apply the standard way (see for instance [6]) to solve Fokker-Planck's conjugate. This requires to solve the Schr\"odinger equation of formula (4.4) below; we saw in section 2 above that its final condition, which is $\phi(s,\cdot)w(s,\cdot)$ for some $\phi\in\tc$, must be in the domain of $\D_\ec$. We shall see that this follows if $\G(w(s,\cdot),w(s,\cdot))$ is bounded; in turn, on $RCD(K,\infty)$ spaces this follows from the fact that $w(s,\cdot)$ is Lipschitz, which we know from section 2. 

\lem{4.2} Let $a,b\in\dc(\D_\ec)\cap L^\infty(S,m)$ and let us suppose that 
$\G(b,b)\in  L^\infty(S,m)$. Then, $ab\in\dc(\D_\ec)$ and 
$$\D_\ec(ab)=
\D_\ec a\cdot b+\D_\ec b\cdot a+2\G(a,b)\in L^2(S,m)  .  \eqno (4.1)$$

\proof Let $a$, $b$ be as above and let $g\in\dc(\ec)\cap L^\infty(S,m)$; from the Leibnitz rule (1.7) we get that
$$\int_S\G(ab,g)\dr m=
\int_S[
\G(a,gb)+\G(b,ga)-2\G(a,b)g
]  \dr m   .   $$
Using (1.8), we get that 
$$\int_S\G(ab,g)\dr m=
\int_S[
-\D_\ec a\cdot b-\D_\ec b\cdot a-2\G(a,b)
]g\dr m     .    \eqno (4.2)$$
This formula holds for $g\in\dc(\ec)\cap L^\infty$, which is dense in $\dc(\ec)$ for the graph norm of $\dc(\ec)$ (with some overkill, this follows by proposition 4.10 of [4]). Thus, by (1.8), $(ab)\in\dc(\D_\ec)$ and (4.1) holds if we show that  
$\D_\ec a\cdot b+\D_\ec b\cdot a+2\G(a,b)\in L^2(S,m)$. 

Now $\D_\ec a\cdot b$ and $\D_\ec b\cdot a$ belong to $L^2$ since 
$\D_\ec a,\D_\ec b\in L^2$ and $a,b\in L^\infty$; moreover, 
$$|\G(a,b)|\le\G(a,a)^\2\G(b,b)^\2  $$
by Cauchy-Schwarz. Since $\G(a,a)\in L^1$ and $\G(b,b)\in L^\infty$ by our hypotheses on $a$ and $b$, we are done. 

\fin

Let $w$ solve equation (4) of the introduction and let $G\in\dc(\D_\ec)\cap L^\infty(S,m)$; we recall that $w(s,\cdot)\in\dc(\D_\ec)$ by lemma 2.1, is bounded by corollary 2.3 and is Lipschitz by lemma 2.4; we can apply lemma 4.2 and get that 
$G\cdot w(s,\cdot)\in\dc(\D_\ec)\cap L^\infty(S,m)$. Applying lemma 2.1 to the final condition $F\cdot w(s,\cdot)$, we find 
$$\psi^s\in C^1((-\infty,s],L^2(S,m))\cap C((-\infty,s],\dc(\D_\ec))     \eqno (4.3)$$
which solves the Schr\"odinger equation, backward in time, 
$$\left\{
\eqalign{
\partial_t\psi^s+&\2\D_\ec\psi^s+\fc\psi^s=0,\quad t\le s\cr
\psi^s(s,\cdot)&= G\cdot w(s,\cdot)  .
}
\right.    \eqno (4.4)$$
By (2.14) we can define $f^s(G,t,x)=\frac{\psi^s(t,x)}{w(t,x)}$; now 
$\frac{1}{w(t,\cdot)}\in\dc(\D_\ec)$ by (1.9) and lemma 2.4; $f^s(G,t,\cdot)\in\dc(\D_\ec)$ by (4.3) and lemma 4.2. More precisely, we have that 
$$f^s(G,\cdot)\in C^1((-\infty,0],L^2(S,m))  \cap
L^\infty_{loc}((-\infty,0],\dc(\D_\ec))   .  $$
Moreover, lemma 2.2 yields the representation formula 
$$f^s(G,t,x)=
\frac{
E^{(t,x)}\left[
e^{
\int_t^s\fc(\tau,\g_\tau)\dr\tau
}    G(\g_s)w(s,\g_s)
\right]
}{
E^{(t,x)}\left[
e^{
\int_t^s\fc(\tau,\g_\tau)\dr\tau
}w(s,\g_s)
\right]
}   =   $$
$$\frac{1}{
w(t,x)
}    \cdot
\int_SG(y)w(s,y)p_{|t-s|}(x,\dr y)
\cdot
\tilde E^{(s-t,x,y)}\left[
e^{
\int_t^s\fc(\tau,\g_\tau)\dr\tau
}
\right]   \eqno (4.5)$$
where we have denoted by $\tilde E^{(s-t,x,y)}$ the expectation with respect to the Brownian bridge which is at $x$ at time $t$ and at $y$ at time $s$; the transition probability 
$p_h(x,\dr y)$ has been defined in section 1. Note that the Brownian bridge exists whenever the Brownian motion exists: it suffices to disintegrate Wiener's measure with respect to the evaluation map $e_0$ of section 1. 

Note that (4.5) defines a function $f^s$ also when $G$ is only continuous: we shall use this fact to define the curve of probability measures. 

We list below the properties of $f^s(G,\cdot)$.

\lem{4.3} Let $(S,d,m)$ satisfy the hypotheses of the $RCD(K,\infty)$ setting and let (F) and (FF) hold. Then the following two points hold. 

\noindent 1) Let $G\in\dc(\D_\ec)\cap L^\infty(S,m)$; then, for all $t\le s$, we have that 
$$f^s(G,\cdot)\in C^1((-\infty,0],L^2)\cap L^\infty_{loc}((-\infty,0],\dc(\D_\ec))$$ 
and 
$$\partial_t f^s(G,t,x)+
\2\D_\ec f^s(G,t,x)+\G(f^s(G,t,\cdot),\log w(t,\cdot))(x)=0
\txt{$m$ a. e. in $S$.}  \eqno (4.6)  $$
Moreover, $f^s(s,x)=G(x)$ for $m$ a. e. $x\in S$. 

\noindent 2) There is a bounded, increasing function 
$\fun{D_6}{[0,+\infty)}{[0,+\infty)}$ such that  for all continuous $G$ we have the following. Let $-T\le t\le s\le 0$; then,  
$$||f^s(G,t,\cdot)||_{L^1(S,m)}\le D_6(T)||G||_{L^1(S,m)}  .  \eqno (4.7)$$

\noindent 3) If $G$ is continuous, then $f^s$ is a semigroup in the past; in other words, if $t\le s\le r$, then
$$f^r(G,t,x)=f^s(
f^r(G,s,\cdot),t,x
)   \txt{for $m$ a. e. $x\in S$.}  \eqno (4.8)$$

\noindent 4) Let $G\in\dc(\D_\ec)\cap L^\infty(S,m)$; then we have that $f^s(G,t,\cdot)$ is differentiable in $s$ and $t$ as a $L^2$-valued function and 
$$\partial_s f^s(G,t,x)|_{t=s}=
-\partial_t|_{t=s} f^s(G,t,x)
\txt{for $m$ a. e. $x\in S$.}    \eqno (4.9)$$

\proof First of all, we recall that $\psi^s\in\dc(\D_\ec)$ by lemma 2.1 while 
$\frac{1}{w}\in\dc(\D_\ec)$ by (1.9), lemma 2.1, corollary 2.3 and lemma 2.4; both functions are bounded by corollary 2.3, and $\frac{1}{w}$ is Lipschitz. Thus, we can apply (4.1) to $a=\psi^s$ and $b=\frac{1}{w}$ and get that 
$$\D_\ec f^s=
\frac{\D_\ec\psi^s}{w}-
2\frac{\G(w,\psi^s)}{w^2}-
\frac{\psi^s}{w^2}\D_\ec w+
\frac{2\psi^s}{w^3}\G(w,w)  .  $$
Now (4.6) follows from the formula above, (4) and (4.4). 

Next, to the bounds on the density of $f^s(\psi,t,\cdot)$. 

Since $\fc$ is bounded, (4.5) and (2.14) imply that there is an increasing function 
$\fun{D_7}{[0,+\infty)}{[0,+\infty)}$ such that, for $t\in[-T,0]$, we have the formula below. 
$$|f^s(G,t,x)|\le D_7(T)\int_S |G(y)| \check p_{|t-s|}(x,y)\dr m(y)  .  $$
Together with Fubini, this yields (4.7) of point 2).  

We only sketch the standard proof of point 3). We start from the right hand side of (4.8); the first and second equalities are (4.5) for $f^s(\cdot,t,\cdot)$ and $f^r(\cdot,s,\cdot)$ respectively; the third equality comes from the fact that the Brownian motion on $[s,r]$ is independent of $[t,s]$; the last equality is (4.5) for $f^r$.
$$f^s(f^r(G,s,\cdot),t,x)=
\frac{
E^{(t,x)}\left[
e^{\int_t^s\fc(\tau,\g(\tau))}f^r(G,s,\g_s)w(s,\g_s)
\right]
}{
w(t,x)
}   =   $$
$$\frac{
E^{(t,x)}\left\{
e^{\int_t^s\fc(\tau,\g(\tau))}    E^{(s,x)}
\left[
e^{\int_s^r\fc(\tau,\g_\tau)\dr\tau}G(\g_r)w(r,\g_r)
\right]
\right\}
}{
w(t,x)
}    =   $$
$$\frac{
E^{(t,x)}
\left[
e^{\int_t^r\fc(\tau,\g_\tau)}G(\g_r)w(r,\g_r)
\right]
}{
w(t,x)
}    =   
f^r(G,t,x)   .   $$

As for point 4), we see from lemma 2.1, corollary 2.3, lemmas 2.4 and 4.2 that the $L^2$-valued function $f^s(G,s,\cdot)$ is differentiable both in $t$ and in $s$. Thus, (4.9) follows from the following formula.
$$0=\frac{\dr}{\dr s}G=
\frac{\dr}{\dr s} f^s(F,s,\cdot)=
\partial_\tau|_{\tau=s}f^s(G,\tau,\cdot)+\partial_\sigma|_{\sigma=t}f^\sigma(G,s,\cdot)  .  $$

\fin

\noindent{\bf End of the proof of proposition 4.1.} {\bf Step 1.} We define the curve of measures.

Let $\fun{G}{S}{\R}$ be continuous; for $t\le s\le 0$, let $f^s$ be defined as in (4.5); we saw above that $f^s(G,t,x)$ is defined for all $s>t$ and all $x\in S$. We can define 
$$\fun{\Lambda_s^{(t,x)}}{C(S)}{\R},\qquad
\fun{\Lambda_s^{(t,x)}}{
G
}{
f^s(G,t,x)
}   .    $$
It is immediate from (4.5) that $\Lambda_s^{(t,x)}$ is linear and brings the non-negative cone into itself; moreover, by lemma 2.2 $\Lambda_s^{(t,x)}1=1$ for all $t\le s\le 0$. Thus, 
$\Lambda_s^{(t,x)}$ is a bounded, positive operator on $C(S)$; arguing as in lemma 2.4 we see that it depends continuously on $s$, $t$ and $x$. Thus, by the Riesz representation theorem, there are measures $\mu^{(t,x)}_s$ depending continuously on $t\le s$ and $x\in S$ such that the second equality below holds; the first one is the definition of $\Lambda^{(t,x)}_s$. 
$$f^s(G,t,x)=
\Lambda_s^{(t,x)}G=\int_S G(y)\dr\mu^{(t,x)}_s(y)  .  \eqno (4.10)$$
Given an initial probability density $\r_t$, we define $\mu^t_s$ as 
$$\mu^t_s=
\int_S\mu_s^{(t,x)}\r_t(x)\dr m(x)  .  \eqno (4.11)$$

\noindent{\bf Step 2.} We assert that $\mu^{(t,x)}_s$ has the semigroup property, forward in time; in other words, if $t<s<r<0$, then
$$\mu^{(t,x)}_{r}=H^s_r\mu^{(t,x)}_s\colon=
\int_S\mu_r^{(s,y)}\dr\mu^{(t,x)}_{s}(y)  .  \eqno (4.12)$$
We briefly prove this fact. Let $G\in C(S,\R)$; the first equality below is (4.10); the second one comes from the fact that the solution of (4.6) is a semigroup in the past, i. e. (4.8). The third equality is (4.10) applied to $f^r$ and $f^s$. 
$$\int_SG\dr\mu^{(t,x)}_{r}=
f^{r}(G,t,x)=$$
$$f^s(
f^r(G,s,\cdot),t,x
)   =
\int_S\dr\mu_s^{(t,x)}(y)
\int_SG(z)\dr\mu_r^{(s,y)}(z)\dr z    .   $$

\noindent{\bf Step 3.} We prove point 1) of proposition 4.1, i. e. that $\fun{}{s}{\mu^t_s}$ satisfies (3.3). 

It suffices to prove that there is $M>0$ such that, for all continuous probability densities 
$G_0$, we have that
$$\left\vert
\int_SG_0(x)\dr\mu_t^s(x)
\right\vert   \le M   .   \eqno (4.13)$$
By (4.11) we get the first equality below, while the second one is (4.10).
$$\int_SG_0(x)\dr\mu^t_s(x)=
\int_S\r_t(x)\dr m(x)\int_SG_0(z)\dr\mu_s^{(t,x)}(z)=$$
$$\int_S f^s(G_0,t,x)\r_t(x)\dr m(x)  .  $$
Now (4.13) follows from (4.7) and the fact that $\r_t$  is bounded. 

\noindent{\bf Step 4.} We show point 2) of proposition  4.1, i. e. that $\fun{}{s}{\mu^{t}_s}$ satisfies (3.4) for $s\in(t,0)$. Let 
$\phi\in\tc$; for simplicity, we shall suppose that $\phi$ is compactly supported in $(t,0)$. For $s\in(t,0)$ we set $G=\phi(s,\cdot)$. For this $G$ we define $f^s$ as in lemma 4.2. The first equality below comes from the fact that $\phi\in\tc$ and dominated convergence.
$$\int_t^0\dr s\int_S  \partial_s\phi(s,z)\dr\mu^{t}_s(z)=$$
$$\lim_{h\searrow 0}
\int_t^0\dr s\int_S  \frac{\phi(s,z)-\phi(s-h,z)}{h}\dr\mu^{t}_s(z)=$$
$$-\lim_{h\searrow 0}
\int_t^0\dr s\int_S\phi(s,z)\dr\left[
\frac{\mu^{t}_{s+h}-\mu^{t}_{s}}{h}
\right]  (z)=  \eqno\hbox{by (4.11)}$$
$$-\lim_{h\searrow 0}\int_t^0\dr s\int_S\r_t(y)\dr m(y)
\int_S\phi(s,z)\dr\left[
\frac{\mu_{s+h}^{(t,y)}-\mu_{s}^{(t,y)}}{h}
\right]  (z)=\eqno\hbox{by (4.12)}$$
$$-\lim_{h\searrow 0}\int_t^0\dr s\int_S\r_t(y)\dr m(y)
\int_S\phi(s,z)\dr\left[
\frac{H^s_{s+h}\mu_{s}^{(t,y)}-\mu_{s}^{(t,y)}}{h}
\right]  (z)=\eqno\hbox{by (4.10)}$$
$$-\lim_{h\searrow 0}
\int_t^0\dr s
\int_S\r_t(y)\dr m(y)
\int_S\frac{f^{s+h}(\phi(s,\cdot),s,z)-\phi(s,z)}{h}
\dr\mu^{(t,y)}_s(z)=\eqno\hbox{by (4.9)}$$
$$\int_t^0\dr s
\int_S\r_t(y)\dr m(y)
\int_S
\partial_t|_{t=s}f^s(\phi(s,\cdot),t,z)\dr\mu^{(t,y)}_s(z)
\eqno\hbox{by (4.6)}$$
$$-\int_t^0\dr s\int_S\r_t(y)\dr m(y)\int_S   [
\2\D_\ec f^s(\phi(s,\cdot),s,x)+\G(f^s(\phi(s,\cdot),s,\cdot),-u(t,\cdot))(x)
]    \dr\mu_s^{(t,x)}(y)=  \eqno{\hbox{by (4.11)}}$$
$$-\int_t^0\dr s
\int_S[
\2\D_\ec\phi(s,y)+\G(\phi(s,\cdot),-u(t,\cdot))(y)
]   \dr\mu_s^{t}(y)  .  $$
Since this is (3.3), we are done. 

\fin

\vskip 2pc
\centerline{\bf Bibliography}

%\noindent [1] R. Adams, J. J. F. Fournier, Sobolev spaces, Academic Press, Singapore, 2009. 

%\noindent [2] A. Ambrosetti, G. Prodi, A primer of nonlinear analysis, Cambridge University Press, Cambridge, 1995.

%\noindent [1] L. Ambrosio, W. Gangbo, Hamiltonian ODE's in the Wasserstein space of probability measures, Communications on Pure and Applied Math., {\bf 61}, 18-53, 2008.

%\noindent [3] L. Ambrosio, J. Feng, On a class of first order Hamilton-Jacobi equations in metric space, preprint. 

\noindent [1] L. Ambrosio, N. Gigli, G. Savar\'e, Gradient Flows, Birkh\"auser, Basel, 2005.

\noindent [2] L. Ambrosio, N. Gigli, G. Savar\'e, Heat flow and calculus on metric measure spaces with Ricci curvature bounded below - the compact case, Analysis and numerics of Partial Differential Equations, 63-115, Springer, Milano, 2013.  

\noindent [3] L. Ambrosio, N. Gigli, G. Savar\'e, Calculus and heat flows in metric measure spaces and applications to spaces with Ricci bounds from below, Invent. Math., {\bf 195}, 289-391, 2014.

\noindent [4] L. Ambrosio, N. Gigli, G. Savar\'e, Metric measure spaces with Riemannian Ricci curvature bounded from below, Duke Math. J., {\bf 163}, 1405-1490, 2014. 

\noindent [5] L. Ambrosio, N. Gigli, G. Savar\'e, Bakry-\'Emery curvature-dimension condition and Riemannian Ricci curvature bounds, Ann. Probab., {\bf 43}, 339-404, 2015. 

\noindent [6] N. Anantharaman, On the zero-temperature or vanishing viscosity limit for certain Markov processes arising from Lagrangian dynamics, J. Eur. Math. Soc. (JEMS), 
{\bf 6}, 207-276, 2004. 

\noindent [7] M. T. Barlow, R. F. Bass, The construction of Brownian motion on the Sierpinski carpet, Ann. IHP,  {\bf 25}, 225-257, 1989.

\noindent [8] M. T. Barlow, E. A. Perkins, Brownian motion on the Sierpiski gasket, Probab. Th. Rel. Fields, {\bf 79}, 543-623, 1988. 

\noindent [9] N. Bouleau, F. Hirsch, Dirichlet forms and analysis on Wiener spaces, Berlin, 1991.

\noindent [10] H. Brezis, Analisi Funzionale, Liguori, Napoli, 1986.

\noindent [11] G. Da Prato, Introduction to Stochastic Differential Equations, SNS, Pisa,  1995. 

%\noindent [12] S. N. Ethier, T. G. Kurtz, Markov Processes, Wiley, Hoboken, 2005. 

\noindent [12] J. Feng, T. Nguyen, Hamilton-Jacobi equations in space of measures associated with a system of conservation laws, Journal de Math\'ematiques pures et Appliqu\'ees, {\bf 97}, 318-390, 2012.

\noindent [13] W. H. Fleming, The Cauchy problem for a Nonlinear first order Partial Differential Equation, JDE, {\bf 5}, 515-530, 1969. 

\noindent [14] M. Fukushima, Y. Oshima, M. Takeda, Dirichlet forms and symmetric Markov processes, De Gruyter, G\"ottingen, 2011. 

%\noindent [GV] D. Gomes, E. Valdinoci, Entropy penalization method for Hamilton-Jacobi equations, Advances in Mathematics, {\bf 215}, 94-152, 2007.

\noindent [15] R. Jordan, D. Kinderleher, F. Otto, The variational formulation of the Fokker-Planck equation, SIAM Journal on Mathematical Analysis, {\bf 29}, 1-17, 1998. 

\noindent [16] T. Kato, Perturbation theory for linear operators, Springer, Berlin, 1980. 

\noindent [17] U. Mosco, Composite media and asymptotic Dirichlet forms, J. Functional Analysis, {\bf 123}, 368-421, 1994.

\noindent [18] A. Pazy, Semigroups of linear operators and applications to Partial Differential Equations, Springer, New York, 1983.

\noindent [19] W. Rudin, Real and Complex Analysis, New Delhi, 1983.

\noindent [20] K.-T. Sturm, Metric measure spaces with variable Ricci bounds and couplings of Brownian motions, Festschrift Masatoshi Fukushima, 553-575, Hackensack, NJ, 2015.

%\noindent [Tay1] M. E. Taylor, Partial Differential Equations, Basic Theory, Springer, New York, 1996.

%\noindent [Tay2] M. E. Taylor, Partial Differential Equation III, Nonlinear Equations, Springer, New York, 1996.

\noindent [21] C. Villani, Topics in optimal transportation, Providence, R. I., 2003.

\end

We have proven (4.9) under the additional hypothesis that 
$f\in L^\infty\cap\dc(\D_\ec)$. Since the Dirichlet form is Markovian, $L^\infty\cap\dc(\D_\ec)$ is dense in $\dc(\D_\ec)$. Thus, if $f\in\dc(\D_\ec)$, we can find 
$f_n\in L^\infty\cap\dc(\D_\ec)$ such that $f_n\tends f$ and $\D_\ec f_n\tends \D_\ec f$ in 
$L^2$; in particular,
$$\int_S[
-\G(f_n,g)-2\G(f_n,V)g
]  \dr m  \tends 
\int_S[
-\G(f,g)-2\G(f,V)g
]  \dr m    .  $$
On the other side, $A$ is closed because it is the generator of a contraction semigroup (we called it $P^V_{s,t}$ in section 2); thus,
$$\int_S(Af_n)g\dr m\tends \int_S(Af)g\dr m  .   $$
Since (4.9) holds for $f_n$, the last two formulas imply that it holds for $f$ too.